\documentclass[nospthms]{svjour2m}                    
\usepackage{mathptmx}      
\usepackage{amsfonts}
\usepackage{amsmath}
\usepackage{amssymb}
\usepackage{amsthm}
\usepackage[usenames]{color}

\newtheorem{theorem}{Theorem}
\newtheorem{lemma}[theorem]{Lemma}

\newtheorem{remark}[theorem]{Remark}

\newcommand{\mc}[1]{{\mathcal #1}}

\newcommand{\bb}[1]{{\mathbb #1}}

\theoremstyle{definition}

\newcommand\N{{\mathbb N}}

\newcommand\cX{{\mathcal X}}

\newcommand\fC{{\textsf{C}}}
\newcommand\fO{{\textsf{O}}}

\journalname{JSP}
\begin{document}

\title{Condensation for a fixed number\\ of independent random variables
}


\author{Pablo A.~Ferrari \and Claudio Landim \and Valentin~V.~Sisko}

  
\institute{Pablo A.~Ferrari  \at
IME USP, Caixa Postal 66281, 05315-970 - S\~{a}o Paulo, BRAZIL \\
\email{pablo@ime.usp.br} 
\and
Claudio Landim  \at
IMPA, Estrada Dona Castorina 110, CEP 22460-320 Rio de Janeiro, Brasil\\
CNRS UMR 6085, Universit\'e de Rouen, UMR 6085, Avenue de
l'Universit\'e, BP.12, Technop\^ole du Madrillet, F76801
Saint-\'Etienne-du-Rouvray, France  \\
\email{landim@impa.br} 
\and
Valentin V.~Sisko \at
IMPA, Estrada Dona Castorina 110, CEP 22460-320 Rio de Janeiro, Brasil\\
\email{valentin@impa.br}
}


\maketitle

\begin{abstract}
 A family of $m$ independent identically distributed random variables indexed by
a chemical potential $\varphi\in[0,\gamma]$ represents piles of particles. As
$\varphi$ increases to $\gamma$, the mean number of particles per site converges
to a maximal density $\rho_c<\infty$. The distribution of particles conditioned
on the total number of particles equal to $n$ does not depend on $\varphi$
(canonical ensemble). 
For fixed $m$, as $n$ goes to infinity the canonical ensemble measure
behave as follows: removing the site with the maximal number of particles, the
distribution of particles in the remaining sites converges to the grand
canonical measure with density $\rho_c$; the remaining particles concentrate
(condensate) on a single site.
\end{abstract}

\section{Introduction}

Condensation phenomena appears in many physical systems. 
From a physical point
of view, in Bose-Einstein condensation a large fraction of the atoms
collapses
into the lowest quantum state, which is possible to observe macroscopically.
From a mathematical point of view it can be seen as the spontaneous migration of
a macroscopic number of particles to some region of the space.  
A physical
account of different models and discussion may be found in the recent
articles:
Evans \cite{E2000}, Jeon and March \cite{JM2000}, Jeon, March and Pittel
\cite{JMP2000}, Godr\`eche \cite{G2003}, Grosskinsky, Sch\"utz and Spohn
\cite{GSS2003}, Evans, Majumdar and Zia \cite{EMSZ2004,EMSZ2006}, Evans and
Hanney \cite{EH2005}, Majumdar, Evans and Zia \cite{MEZ2005} and Godr\`eche and
Luck \cite{GL2005}.

These references consider the following mathematical model. Fix a natural number
$m$. The state space consists of particle configurations $\xi\in\N^m$; $\xi_i$
represents the number of particles at site $i\in\{1,\dots,m\}$. 
Consider a family $\mu_\varphi$ of product measures on $\N^m$ indexed
by a chemical potential $\varphi$ with range in $[0,\gamma]$ for some
$\gamma>0$.  The density $\rho = R(\varphi)$ (mean number of particles
per site) is an increasing function of $\varphi$. We assume that the
maximal density is finite: $\rho_c = R(\gamma)<\infty$.

For an integer $n\ge 0$, denote by $\nu_n$ the canonical measure
associated to the family $\mu_\varphi$. 
This is the measure
$\mu_\varphi$ conditioned on the hyperplane of configurations with $n$
particles: $\nu_n (\xi) = \mu_\varphi(\xi \mid \sum_i \xi_i = n)$.
In this context, condensation means that all but a few particles
concentrate on the same site. We prove in this article such a
statement under some conditions on the product measure $\mu_\varphi$.
More precisely, we show that if we remove the site with the largest
number of particles, the projection of the canonical measure on the
remaining sites converges to the product measure with maximal density
$\rho_c$, when the total number of sites $m$ is
fixed and the total number $n$ of particles increases to infinity.

\section{Notation and result}
\label{s:mtp}

Denote by $\bb N$ the nonnegative integers and let $f : \bb N \to \bb
R_+$ be a positive real function such that
\begin{eqnarray}
\label{a2}
\lim_{n\to\infty} \frac{f(n)}{f(n+1)} \;=\; \gamma \;, \quad
\sum_{n\ge 0} \gamma^n f(n) \;<\; \infty
\end{eqnarray}
for some $\gamma >0$. Assume also that $\gamma ^n f(n)$ is decreasing
and that there exists a function $C : \bb N\to\bb R$ such that
\begin{equation}
\label{a1}
\gamma ^{n/m} f(n/m) \;\le\; C (m) \gamma^{n} f(n)
\end{equation}
for all $n$. $f(n) = n^{-\alpha}$ with $\alpha > 1$, for instance,
satisfies \eqref{a2}, \eqref{a1}. 

We consider a family of product measures in the state space $\bb N^m$.
Configurations in this space are denoted by the Greek letter $\xi = (\xi_1 ,
\dots , \xi_m)$, with $\xi_j$ in $\bb N$ for $1\le j\le m$.  Let $Z:\bb R_+ \to
\bb R_+$ the partition function defined by
\begin{equation*}
Z(\varphi) \;=\; \sum_{n\ge 0} \varphi^n f(n)\;.
\end{equation*}
It follows from \eqref{a2} that the radius of convergence of $Z$ is
equal to $\gamma$: $Z(\varphi)<\infty$ if and only if $\varphi\le
\gamma$. For $0\le \varphi\le \gamma$, denote by $\mu_{m,\varphi}$ the
grand canonical measure on $\bb N^m$ given by
\begin{equation*}
\mu_{m,\varphi} (\xi) \;=\; \frac{1}{Z(\varphi)^m}
\prod_{j=1}^m \varphi^{\xi_j} f(\xi_j)\;.
\end{equation*}

Denote by $|\xi|$ the total number of particles of the configuration:
$|\xi| = \xi_1 + \cdots + \xi_m$ and by $\Sigma_{m,n}$, $n\ge 1$, the
subspace of configurations of $\bb N^m$ with $n$ particles:
\begin{equation}
\label{a3}
\Sigma_{m,n} \;=\; \{\xi\in\N^m \,:\, |\xi|=n\}\;.
\end{equation}
For $n\ge 0$, let $\nu_{m,n}$ be the canonical measure concentrated on
configurations with $n$ particles:
\begin{equation*}
\nu_{m,n}(\xi)\;=\; \mu_{m,\varphi} (\xi\,|\,|\xi|=n) \;.
\end{equation*}
Notice that the right hand side does not depend on the parameter
$\varphi$. 

Denote by $\mc X_m$ the space of monotone configurations:
\begin{equation*}
\cX_m \;=\; \{\zeta \in \N^m:\zeta(1)\le\dots\le\zeta(m)\}\;.
\end{equation*}
Configurations of $\mc X_m$ are denoted by the Greek letter $\zeta$.  Define the
order operator $\fO : \bb N^m \to \mc X_m$ which takes a configuration of $\bb
N^m$ to the monotone configuration by reordering appropriately the coordinates:
\begin{eqnarray*}
\!\!\!\!\!\!\!\!\!\!\!\!\!\!\! &&
(\fO \xi)_i \;\le\; (\fO \xi)_j \quad \text{for $1\le i\le j\le m$} \\
\!\!\!\!\!\!\!\!\!\!\!\!\!\!\! && \quad 
\text{and } (\fO \xi)_i = \xi_{\sigma(i)} \quad\text{for 
a permutation $\sigma$ of the set $\{1,\dots ,m\}$} \; .
\end{eqnarray*}
Define the cut operator $\fC:\bb N^m \to \bb N^{m-1}$ which eliminates
the last coordinate:
\[
(\fC \xi)_i  = \xi_i \quad \text{for $1\le i\le m-1$.}
\]

These operators induce two measures: $\hat \mu_{m, \varphi}$ defined on
$\cX_{m}$ and $\hat\nu_{m-1,n}$ defined on $\cX_{m-1}$:
\begin{equation*}
\hat \mu_{m,\varphi} \;=\; \mu_{m,\varphi}  \fO^{-1} \; , \qquad  
\hat \nu_{m-1,n} = \nu_{m,n} (\fC\circ\fO)^{-1} \;.
\end{equation*}
To obtain a sample of $\hat \mu_{m, \varphi}$, sample a configuration with the
product measure $\mu_{m,\varphi}$ and order it.  To obtain a sample of
$\hat\nu_{m-1,n}$, sample a configuration with the product measure
$\mu_{m,\varphi}$ conditioned to have $n$ particles (i.e.  $\nu_{m,n}$), order
it and cut the last coordinate.

Next we state our result. It asserts that if the total number of sites $m$ is
fixed and the total number $n$ of particles increases to infinity, the canonical
measure concentrated on configurations with $n$ particles, from which we removed
the site with the largest number of particles, converges to the ordered
\emph{maximal} grand canonical measure. Maximal means that we set $\varphi$ to
be $\gamma$.

By symmetry the position of the maximum site is uniformly distributed in
$\{1,\dots,m\}$. Sampling from the ordered measure $\hat\nu_{m-1,n}$ and
reordering randomly and uniformly the labels of the particles we obtain a
measure which approachs the maximal product measure.

\begin{theorem}
\label{t1}
Assume that f is a nonincreasing function satisfying hypotheses \eqref{a2} and
\eqref{a1}. For each $m\ge 2$, the sequence of probability measures
$\hat\nu_{m-1,n}$ converges weakly to $\hat \mu_{m-1, \gamma}$ as $n$ tends to
infinity.
\end{theorem}

\begin{remark}
  \label{t2} {\rm A similar result has been proven in various of the quoted
    references in the thermodynamical limit. See for instance
    \cite{GSS2003}, where the the number of particles $n$ divided the
    number of sites $m=m(n)$ converges to a constant $\rho>\rho_c$ as
    $n$ goes to infinity.  Theorem \ref{t1} shows that this phenomenon
    is just a combinatorial fact that can be observed already without
    making the number of sites grow to infinity.  }
\end{remark}

\begin{remark}
\label{t3}
{\rm The condition that $f$ is decreasing can be relaxed. We can replace
this condition and \eqref{a1} by the hypothesis that for every $m\ge 1$,
there exists a finite constant $C_1 (m)$ such that 
\begin{equation}
\label{a6}
\frac{\max_{n/m \le k \le n} \gamma ^k f (k)}
{\min_{n/m \le k \le n} \gamma ^k f (k)} \;\le\; C_1(m)
\end{equation}
for all $n\ge 1$.}
\end{remark}

\begin{remark}
  \label{t4} 
  
  {\rm Relation with the zero range process. Fix $m\ge 1$, a
    irreducible symmetric transition probability $p(x,y)$ on $\{1,
    \dots, m\}$ and a positive function $g: \bb N\to\bb R_+$. The
    zero-range process associated to $(p,g)$ can be informally
    described as follows. Particles evolve on $\{1, \dots, m\}$. If
    there are $k$ particles at some site $x$ of $\{1, \dots, m\}$ one
    of them jumps to site $y$ at rate $g(k) p(x,y)$, independently
    from what happens at the other sites.
    
    Let $f(k)^{-1} = g(1) \cdots g(k)$ and assume that $f$ satisfies
    assumption \eqref{a2} for some $\gamma >0$.  It is well known
    \cite{KL1999} that the product measures $\{\mu_{m,\varphi} : 0\le
    \varphi\le \gamma\}$ defined above are invariant for the symmetric
    zero-range process.
    
    Let $R^* = \sum_{n\ge 0} n \gamma^n f(n)$ and assume that
    $R^*<\infty$. The function $R:[0,\gamma] \to [0, R^*]$ defined by
    $R(\varphi) = \mu_{m,\varphi} [\xi_1]$ gives the density of
    particles under the invariant measure $\mu_{m,\varphi}$. A simple
    computation of $R'$ shows that $R$ is bijective. Since in the
    symmetric case any product invariant measure belongs to the set
    $\{\mu_{m,\varphi} : 0\le \varphi\le \gamma\}$, there are no
    product invariant measures with density above $R^*$.

    Since the process is ergodic, for every $n\ge 1$, there exists a unique
    stationary measure, denoted above by $\nu_{m,n}$. Order the configuration
    and remove the site with the largest number of particles.  Under the
    additional assumption \eqref{a6}, Theorem \ref{t1} states that this measure
    converges, as $n\uparrow\infty$, to the ordered maximal product state.  In
    particular, all but a finite number of particles tend to accumulate in one
    site.}
\end{remark}

\section{Proof of Theorem \ref{t1}}

For an ordered configuration $\zeta$ belonging to $\cX_{m}$ (resp.
$\eta$ belonging to
$\cX_{m-1}$), let $K_m(\zeta)$ (resp.  $K_{m-1,n}(\eta)$) be the
number of configurations in $\bb N^m$ whose ordering (resp. ordering
and cutting) gives $\zeta$ (resp. $\eta$):
\begin{eqnarray*}
\!\!\!\!\!\!\!\!\!\!\!\!\!\! &&
K_m(\zeta) \;=\; \sharp\{\xi \in \N^{m}:\fO(\xi)=\zeta\}\; , \\
\!\!\!\!\!\!\!\!\!\!\!\!\!\! && \quad
K_{m-1,n}(\eta) \;=\; \sharp\{\xi \in \Sigma_{m,n}:
\fC\circ\fO(\xi)=\eta\}\;. 
\end{eqnarray*}
Note that the number of elements of an empty set is zero.  Since
$\mu_{m,\varphi}$, $\nu_{m,n}$ are invariant by permutation, for any
$\zeta$ in $\cX_{m}$, $\eta$ in $\cX_{m-1}$,
\begin{eqnarray}
\label{a4}
\!\!\!\!\!\!\!\!\!\!\!\!\!\! &&
\hat \mu_{m, \varphi} (\zeta) \;=\; K_m(\zeta) \, \mu_{m, \varphi}
(\zeta)\; , \\
\!\!\!\!\!\!\!\!\!\!\!\!\!\! &&
\quad \hat \nu_{m-1,n}(\eta) \;=\;  K_{m-1,n}(\eta) \,
\nu_{m,n}(\eta_1, \dots, \eta_{m-1}, n - |\eta|)\;.
\nonumber
\end{eqnarray}

\begin{lemma}
\label{lm:g}
For any  $\eta$ in $\cX_{m-1}$ and any $n\ge 1$, $K_{m-1,n} (\eta) \le  m
K_{m-1} (\eta)$. Moreover, for any $\eta$ in $\cX_{m-1}$,
\begin{equation*}
\lim_{n\to\infty} K_{m-1,n} (\eta) \;=\;  m K_{m-1} (\eta)\;.
\end{equation*}
\end{lemma}
                          
The proof is elementary. A fixed ordered configuration $(\eta_1 ,
\dots , \eta_{m-1})$ belonging to $\cX_{m-1}$ arises from a permutation of the
coordinates $\eta_j$. To compute $K_{m-1,n}$, we add the coordinate
$n - |\eta|$.  This extra coordinate accounts for the factor $m$ which
corresponds to all its possible positions. If $n$ is large this last
coordinate is different from all others and provides $m$ distinct
vectors.
\medskip
                                                                        
We are now in position to prove Theorem \ref{t1}. Fix $\eta$ in
$\cX_{m-1}$. In view of \eqref{a4},
\begin{equation*}
\hat \nu_{m-1,n} (\eta) \;=\; 
\frac{K_{m-1,n} (\eta) f (\eta_1 ) \cdots  f (\eta_{m-1} ) f (n - |\eta|)}
{\sum_{\eta' \in \cX_{m-1}} K_{m-1,n} (\eta') f (\eta'_1 ) \cdots  f
  (\eta'_{m-1} ) f (n - |\eta'|) }\;\cdot
\end{equation*}
By \eqref{a4} and Lemma \ref{lm:g}, to prove Theorem \ref{t1}, we only
need to show that
\begin{eqnarray}
\label{a5}
\!\!\!\!\!\!\!\!\!\!\!\!\!\!\!\!\!\!\! &&
\lim_{n\to\infty} \sum_{\eta' \in \cX_{m-1}} K_{m-1,n} (\eta') 
f (\eta'_1 ) \cdots  f  (\eta'_{m-1} ) \frac{f (n - |\eta'|)}
{f (n - |\eta|)} \\
\!\!\!\!\!\!\!\!\!\!\!\!\!\!\!\!\!\!\! && \quad
=\; \sum_{\eta' \in \cX_{m-1}} m K_{m-1} (\eta') 
\gamma^{|\eta'| - |\eta|} f (\eta'_1 ) \cdots  f  (\eta'_{m-1} )\;.
\nonumber
\end{eqnarray}

Fix a positive constant $M > 1$. By assumption \eqref{a2} and by Lemma
\ref{lm:g}, for configurations such that $|\eta| \le M$ , $\lim_n
K_{m-1,n} (\eta) = m K_{m-1} (\eta)$ and 
\[
\lim_{n \to \infty} f (n - |\eta|)/f (n -
|\eta'|) = \gamma^{|\eta'| - |\eta|}.
\]
In particular, for every $M\ge 1$,
\begin{eqnarray*}
\!\!\!\!\!\!\!\!\!\!\!\!\!\!\!\!\!\!\! &&
\lim_{n\to\infty} 
\sum_{\substack{\eta' \in \cX_{m-1} \\ |\eta'|\le M}} 
K_{m-1,n} (\eta') f (\eta'_1 ) \cdots  f  (\eta'_{m-1} ) 
\frac{f (n - |\eta'|)} {f (n - |\eta|)} \\
\!\!\!\!\!\!\!\!\!\!\!\!\!\!\!\!\!\!\! && \quad
=\; \sum_{\substack{\eta' \in \cX_{m-1} \\ |\eta'|\le M}} 
m K_{m-1} (\eta') \gamma^{|\eta'| - |\eta|} f (\eta'_1 ) 
\cdots  f (\eta'_{m-1} )\;.
\end{eqnarray*}

To estimate the sum $|\eta'| > M$, recall from Lemma \ref{lm:g} that 
\[
K_{m-1,n} (\eta') \le m K_{m-1} (\eta').
\]
On the other hand, since
\[
\eta_m = n - |\eta|\ge \max_{1\le j \le m-1} \eta_j, \quad
\eta'_m = n - |\eta'|\ge \max_{1\le j \le m-1} \eta'_j,
\]
 we have that $n/m\le
\eta_m$ , $\eta'_m\le n$. Thus, by assumption \eqref{a6}, 
\begin{equation*}
  \frac{f (n - |\eta'|)}{f (n - |\eta |)} \;\le\; C_1(m) 
\gamma^{|\eta'| - |\eta|}.
\end{equation*}
Therefore,
\begin{eqnarray*}
\!\!\!\!\!\!\!\!\!\!\!\!\!\!\!\!\!\!\! &&
\sum_{\substack{\eta' \in \cX_{m-1} \\ |\eta'|> M}} 
K_{m-1,n} (\eta') f (\eta'_1 ) \cdots  f  (\eta'_{m-1} ) 
\frac{f (n - |\eta'|)} {f (n - |\eta|)} \\
\!\!\!\!\!\!\!\!\!\!\!\!\!\!\!\!\!\!\! && \quad
\le \; C(m) \sum_{\substack{\eta' \in \cX_{m-1} \\ |\eta'|> M}} 
K_{m-1} (\eta') \gamma^{|\eta'| - |\eta|} f (\eta'_1 ) 
\cdots  f (\eta'_{m-1} )\;.
\end{eqnarray*}
This expression vanishes as $M\uparrow\infty$ because the sum without
the constraint $|\eta'|> M$ is finite, being equal to
$\gamma^{-|\eta|} Z(\gamma)^{m-1}$, where $Z(\gamma)$ is the
partition function introduced in Section \ref{s:mtp}.  This concludes
the proof of the theorem.


\begin{acknowledgements}
The authors are grateful to CNPq, FAPERJ, FAPESP, and PRONEX for financial support. 
\end{acknowledgements}


\end{document}